\documentclass[12pt]{article}
\usepackage{amsfonts,amsmath,amsthm, amssymb}
\usepackage{latexsym, euscript, epic, eepic}
\renewcommand{\theequation}{\thesection.\arabic{equation}}
\newtheorem{cro}{Corollary}[section]
\newtheorem{defn}{Definition}[section]
\newtheorem{prop}{Proposition}[section]
\newtheorem{thm}{Theorem}[section]
\newtheorem{lem}{Lemma}[section]
\newtheorem{rem}{\bf Remark}[section]

\begin{document}
\title{The variational principle of local pressure for actions of
sofic group
 }
\author{Xiaoyao Zhou ,$^\dag$  Ercai Chen$^{\dag \ddag}$\\ 
\small \it $\dag$ School of Mathematical Science, Nanjing Normal University,\\
\small \it Nanjing 210097, Jiangsu, P.R.China\\
\small \it $\ddag$ Center of Nonlinear Science, Nanjing University,\\
\small \it  Nanjing 210093, Jiangsu, P.R.China\\
\small \it  e-mail: zhouxiaoyaodeyouxian@126.com  ecchen@njnu.edu.cn\\
}
\date{}
\maketitle

\begin{center}
\begin{minipage}{120mm}
{\small {\bf Abstract.} This study establishes the variational
principle for local pressure in the sofic context. }
\end{minipage}
\end{center}

\vskip0.5cm
\section{Introduction}
Variational principles are beautiful results in a dynamical system.
Establishing variational principle is an important topic in
dynamical system theory. The first variational principle that
reveals the relationship between topological entropy and
measure-theoretic entropy was obtained by L. Goodwyn \cite{[Goo2]}
and T. Goodman \cite{[Goo1]}. M. Misiurewicz gave a short proof in
\cite{[Mis]}. R. Bowen \cite{[Bow0]} established the variational
principle of entropy for non-compact set in 1973. Y. Pesin and B.
Pitskel studied the variational principle of pressure for
non-compact set.

Romagnoli \cite{[Rom]} introduced two types of measure-theoretic
entropies relative to a finite open cover and proved the variational
principle between local entropy and one type. Later, Glasner and
Weiss \cite{[GlaWei]} proved that if the system is invertible, then
the local variational principle is also true for another type
measure-theoretic entropy.  W. Huang \& Y. Yi \cite{[HuaYi]}
obtained the variational principle of local pressure. W. Huang, X.
Ye and G. Zhang \cite{[HuaYeZha]} established the variational
principle of local entropy for a countable discrete amenable group
action. B. Liang and K. Yan \cite{[LiaYan]} generalized it to local
pressure for sub-additive potentials of amenable group actions.

Recently, L. Bowen \cite{[Bow2]}, \cite{[Bow1]} introduced a notion
of entropy for measure-preserving actions of a countable discrete
sofic group on a standard probability space admitting a generating
finite partition. Just after that, D. Kerr and H. Li \cite{[KerLi1]}
 established variational principle of entropy for
sofic group actions. G. Zhang \cite{[Zha]} generalized it to the
variational principle of local entropy for countable infinite sofic
group actions. N. Chung \cite{[Chu]} generalized it to the
variational principle of pressure for sofic group actions.

Local sofic pressure is  introduced and variational principle is
established in this study.

\section{Preliminaries}
For a set $Z,$ denote by $\mathcal{F}(Z)$ the set of all non-empty
finite subsets of $Z$ and by $|Z|$ its cardinality. For
$d\in\mathbb{N},$ write $[d]$ for the set $\{1,2,\cdots,d\}$ and
Sym$(d)$ for the permutation group of $[d].$ Fix $G$ to be a
countable infinite discrete sofic group with unit $e$ in this pape
(except some special statements), i.e., there exists a sequence
$\{\sigma_i\}_{i\in\mathbb{N}}$ of maps $\sigma_i:G\to {\rm
Sym}(d_i),g\to \sigma_{i,g},$ which is asymptotically multiplicative
and free in the sense that

\begin{equation*}
\lim\limits_{i\to\infty}\frac{1}{d_i}\left|\{a\in
[d_i]:\sigma_{i,st}(a)=\sigma_{i,s}\sigma_{i,t}(a)\}\right|=1
\end{equation*}
for all $s,t\in G$ and

\begin{equation*}
\lim\limits_{i\to\infty}\frac{1}{d_i}\left|\{a\in
[d_i]:\sigma_{i,s}(a)\neq\sigma_{i,t}(a)\}\right|=1
\end{equation*}
for all $s\neq t\in G.$

Such a sequence $\Sigma=\{\sigma_i\}_{i\in \mathbb{N}}$ with
$\lim\limits_{i\to\infty}d_i=\infty$ is referred as a sofic
approximation sequence of $G$ and fix $\Sigma$ in this paper.

Let $(X,\rho)$ be a compact metric space. Denote by $M(X)$ the set
of all Borel probability measures on $X,$ by $M(X,G)$ the set of all
$G-$invariant Borel probability measures on $X,$ by $M^e(X,G)$ the
set of all ergodic measures on $X,$ respectively.

Let $C(X)$ be the set of all real valued continuous functions on
$X.$ Define norm $\|\cdot\|_\infty$ on $C(X)$ by
$\|f\|_\infty=\sup\limits_{x\in X}|f(x)|.$ The actions of $G$ on
points will usually be expressed by the concatenation $(s,x)\to sx.$
For a map $\sigma:G\to {\rm Sym}(d)$ for some $d\in \mathbb{N},$
denote $\sigma_s(a)$ for $s\in G$ and $a\in [d]$ simply by $sa$ when
convenient. A continuous action $\alpha$ of $G$ on $(X,\rho)$
induces an action on $C(X)$ as follows: for $g\in C(X)$ and $s\in
G,$ the function $\alpha_s(g)$ is given by $x\to g(s^{-1}x).$ For
$Y\subset X,$ denote by $C_Y$ the set of covers of $Y,$ by $C_Y^o$
the set of open covers of $Y$(Open sets in $C_Y^o$ are with respect
to $X.$), and by $\mathcal{P}_X$ the set of all finite Borel
partitions of $X.$ Let $\mathcal{U}\in C_X^o$ and $d\in\mathbb{N},$
set $\mathcal{U}^d=\{U_{i_1}\times U_{i_1}\times \cdots\times
U_{i_d}: U_{i_j}\in\mathcal{U},j=1,\cdots,d\}.$ For two given covers
$\mathcal{U}, \mathcal{V}\in C_X,\mathcal{U}$ is said to be finer
than $\mathcal{V}$ (denoted by $\mathcal{U}\succcurlyeq\mathcal{V}$
or $\mathcal{V}\preccurlyeq\mathcal{U}$) if each element of
$\mathcal{U}$ is contained in some element of $\mathcal{V}.$ For
$\mathcal{U}\in C_X$ and $\emptyset\neq Y\subset X,$ we set
$N(\mathcal{U},Y)$ to be the minimal cardinality of sub-families of
$\mathcal{U}$ covering $Y.$ Set $N(\mathcal{U},\emptyset)=0$ by
convention.

Let $F\in \mathcal{F}(G)$ and $\delta>0.$ Let $\sigma$ be a map from
$G$ to Sym$(d)$ for some $d\in\mathbb{N}.$

\begin{defn}{\rm\cite{[KerLi1]}\cite{[KerLi2]}\cite{[Li]}}Define

$ Map^\rho(F,\delta,\sigma) =\left\{\varphi:[d]\to X:
\max\limits_{s\in
F}\rho_2(\varphi\circ\sigma_s,\alpha_s\circ\varphi)<\delta\right\},$
where
\begin{equation*}
\rho_2(\varphi,\psi)=\left(\frac{1}{d}\sum\limits_{a=1}^{d}(\rho(\varphi(a),\psi(a)))^2\right)^{1/2}.
\end{equation*}
\end{defn}

\noindent There is another compatible metric on $ X^d$ by
\begin{equation*}
\rho_\infty(\varphi,\psi)=\max\limits_{a\in
[d]}\rho(\varphi(a),\psi(a)).
\end{equation*}

\begin{lem}\label{lem2.1}{\rm\cite{[Zha]}}
Let $\rho_1$ and $\rho_2$ be two compatible metrics on $X$ and
$F\in\mathcal{F}(G).$ Then for each $\delta_2>0$ there exists
$\delta_2\geq\delta_1>0$ such that

\begin{equation*}
Map^{\rho_1}(F,\delta_1,\sigma)\subset
Map^{\rho_2}(F,\delta_2,\sigma)
\end{equation*}
for each map $\sigma: G\to$ Sym($d$) with some $d\in\mathbb{N}.$
\end{lem}

\begin{defn}{\rm\cite{[KerLi1]}\cite{[KerLi2]}} Let $L\in \mathcal{F}(C(X))$ and $\mu\in M(X).$ Set
\begin{equation*}
 Map^\rho_\mu(L,F,\delta,\sigma)=\left\{\varphi\in
Map^\rho(F,\delta,\sigma):\max\limits_{f\in
L}\left|\frac{1}{d}\sum\limits_{j=1}^d f(\varphi(j))-
\mu(f)\right|<\delta\right\}.
\end{equation*}
\end{defn}

\begin{defn}
Let $\mathcal{U}\in C_X^o$ and $f\in C(X).$ Set

\begin{eqnarray*}\begin{split}
&N^\rho(\sigma,F,\delta,\mathcal{U},f)=\inf
\{\sum\limits_{U\in\beta}\sup\limits_{\varphi\in U\cap
Map^\rho(F,\delta,\sigma)}\exp(\sum\limits_{a=1}^{d}f(\varphi(a))):\beta\in
C_{Map^\rho(F,\delta,\sigma)}{\rm~ and
~}\beta\succcurlyeq\mathcal{U}^d\},
\end{split}\end{eqnarray*}

\begin{equation*}
P^\rho(F,\delta,\mathcal{U},f)=\limsup\limits_{i\to\infty}\frac{1}{d_i}\log
N^\rho(\sigma_i,F,\delta,\mathcal{U},f),
\end{equation*}
\begin{equation*}
P(\mathcal{U},f)=\inf\limits_{F\in
\mathcal{F}(G)}\inf\limits_{\delta>0}P^\rho(F,\delta,\mathcal{U},f),
\end{equation*}
\begin{equation*}
P(X,f)=\sup\limits_{\mathcal{U}\in C_X^o} P(\mathcal{U},f).
\end{equation*}
\end{defn}

\begin{rem}
\noindent {\rm (i)} By Lemma \ref{lem2.1}, we have that
$P(\mathcal{U},f)$ does not depend on the compatible metric on $X.$
We call it local sofic pressure with respect to $\mathcal{U}$ and
$f.$

\noindent {\rm (ii)} Simplify  $Map^\rho(F,\delta,\sigma)$ as
$Map(F,\delta,\sigma),$  write $N(\sigma_i,F,\delta,\mathcal{U},f)$
instead of $N^\rho(\sigma_i,F,\delta,\mathcal{U},f),$ and abridge
$P^\rho(F,\delta,\mathcal{U},f)$ to $P(F,\delta,\mathcal{U},f),$
although these quantities depend on the compatible metric $\rho.$
\end{rem}

\begin{defn}
Let $\mathcal{U}\in C_X^o$ and $f\in C(X).$ Set
\begin{eqnarray*}\begin{split}
&N^*(\sigma,F,\delta,\mathcal{U},f)=\inf
\{\sum\limits_{U\in\beta}\sup\limits_{\varphi\in U\cap
Map(F,\delta,\sigma)}\exp(\sum\limits_{a=1}^{d}f(\varphi(a))):\beta\in
C^o_{Map(F,\delta,\sigma)}{\rm~ and
~}\beta\succcurlyeq\mathcal{U}^d\},
\end{split}\end{eqnarray*}

\begin{equation*}
P^*(F,\delta,\mathcal{U},f)=\limsup\limits_{i\to\infty}\frac{1}{d_i}\log
N^*(\sigma_i,F,\delta,\mathcal{U},f),
\end{equation*}
\begin{equation*}
P^*(\mathcal{U},f)=\inf\limits_{F\in
\mathcal{F}(G)}\inf\limits_{\delta>0}P^*(F,\delta,\mathcal{U},f).
\end{equation*}
\end{defn}

\begin{rem}
\noindent {\rm (i)} $\varphi:[d]\to X$ can be seen as a point
$(\varphi(1),\cdots,\varphi(d))\in X^d$ {\rm (See \cite{[Zha]})}.

\noindent{\rm (ii)} Combining the variational principle in
{\rm\cite{[Chu]}} and the variational principle in this study, we
can obtain that $P(X,f)$ coincides with the sofic pressure in
{\rm\cite{[Chu]}}.

\noindent{\rm (iii)} $P(\mathcal{U},0)= h(G,\mathcal{U}), $ where
the local entropy $h(G,\mathcal{U})$ is defined in
{\rm\cite{[Zha]}}.

\noindent{\rm (iv)} It is obvious that $P^*(\mathcal {U},f)\geq
P(\mathcal {U},f).$ We will prove the opposite inequality by the
proof of the variational principle.
\end{rem}

\begin{defn}{\rm\cite{[Zha]}}
Let $ \mathcal{U}\in C_X^o.$ Define $h_\mu(\mathcal{U})$ to be the
$\mu$-measure-theoretic entropy of $\mathcal{U}$ as follows:
\begin{equation*}
h^\rho_\mu(L,F,\delta,\mathcal{U})=\limsup\limits_{i\to\infty}\frac{1}{d_i}\log
N(\mathcal{U}^{d_i},Map^\rho_{\mu}(L,F,\delta,\sigma_i)),
\end{equation*}

\begin{equation*}
h_\mu(\mathcal{U})=\inf\limits_{L\in
\mathcal{F}(C(X))}\inf\limits_{F\in
\mathcal{F}(G)}\inf\limits_{\delta>0}h^\rho_\mu(L,F,\delta,\mathcal{U}).
\end{equation*}
\end{defn}

\begin{rem}
\noindent {\rm (i)} $h_\mu(\mathcal{U})$ does not depend on the
choice of compatible metric in $X.$

\noindent {\rm (ii)} Abbreviate $Map_\mu(L,F,\delta,\sigma)$ to
$Map^\rho_\mu(L,F,\delta,\sigma),$ and write
$h_\mu(L,F,\delta,\mathcal{U})$ instead of
$h^\rho_\mu(L,F,\delta,\mathcal{U})$ for convenience, though the two
quantities depend on the compatible metric in $X.$
\end{rem}
\section{Variational principle for sofic local pressure}

We give two lemmas to state our main result.
\begin{lem}\label{lem3.1}
If $\mathcal{U}\in C_X^o, f\in C(X),$ then there exists $\mu\in
M(X,G)$ such that $P^*(\mathcal{U},f)\leq
h_\mu(\mathcal{U})+\mu(f).$
\end{lem}
\noindent{\it Proof.} Take a sequence $e\in F_1\subset
F_2\subset\cdots$ with $F_i\in \mathcal{F}(G)$ for all $i$ and
$\bigcup\limits_{i}F_i=G.$ Since $(X,\rho)$ is a compact metric
space, there exists a sequence $\{g_m\}_{m\in\mathbb{N}}$ that
 is dense in $C(X).$ Set $L_n=\{f,g_1,\cdots,g_n\}.$ There exists $Q>0$ such that
$\max\limits_{g\in L_n}\|g\|_{\infty}\leq Q.$

Choose $\delta_n>0$ such that

\noindent (i) $\delta_n<\min\left\{\frac{1}{12Q|F_n|},
\frac{1}{9n}\right\};$

\noindent (ii) $\max\limits_{g\in L_n}|g(x)-g(y)|<\frac{1}{6n}$  for
all $x,y\in X$ with $\rho(x,y)<\sqrt{\delta_n}.$\\

We will find some $\mu_n\in M(X)$ such that

\noindent (i) $h_{\mu_n}(L_n, F_n,
\frac{1}{3n},\mathcal{U})+\mu_n(f)+\frac{10}{9n}\geq
P^*(F_n,\delta_n,\mathcal{U},f);$

\noindent (ii) $\max\limits_{g\in L_n}\max\limits_{t\in
F_n}|\mu_n(\alpha_{t^{-1}}(g))-\mu_n(g)|<\frac{1}{n}.$\\

By weak* compactness there exists a finite subset $D\subset M(X)$
such that for any map $\sigma:G\to {\rm Sym}(d)$ for some $d\in
\mathbb{N}$ and any $\varphi\in Map(F_n,\delta_n,\sigma)$ there is a
$\mu_\varphi\in D$ such that $\max\limits_{g\in
L_n}\max\limits_{t\in
F_n}|\mu_\varphi(\alpha_{t^{-1}}(g))-(\varphi_*\zeta)(\alpha_{t^{-1}}(g))|<\delta_n,$
where
$(\varphi_*\zeta)(h)=\frac{1}{d}\sum\limits_{a=1}^{d}h(\varphi(a))$
for all $h\in C(X).$

 For each $\varphi\in Map(F_n, \delta_n, \sigma),$
set $\Lambda_\varphi=\left\{a\in [d]:\max\limits_{t\in
F_n}\rho(\varphi(ta),t\varphi(a))<\sqrt{\delta_n}\right\}.$ Then
$|\Lambda_\varphi|\geq(1-|F_n|\delta_n)d.$

Thus,
\begin{equation*}\begin{split}
&\max\limits_{g\in L_n}\max\limits_{t\in
F_n}|(\varphi_*\zeta)(\alpha_{t^{-1}}(g))-((\varphi\circ
\sigma_t)_*\zeta)(g)|\\&\leq \max\limits_{g\in L_n}\max\limits_{t\in
F_n}\frac{1}{d}\left|\sum\limits_{a\in
\Lambda_\varphi}(g(t\varphi(a))-g(\varphi(ta)))\right|+\frac{1}{d}\left|\sum\limits_{a\notin
\Lambda_\varphi}(g(t\varphi(a))-g(\varphi(ta)))\right|\\&\leq
\frac{1}{d}|\Lambda_\varphi|\cdot\frac{1}{6n}+\frac{1}{d}2Q|F_n|\delta_nd\\&\leq
\frac{1}{6n}+\frac{1}{6n}=\frac{1}{3n},
\end{split}\end{equation*}

and hence
\begin{equation*}\begin{split}
&|\mu_\varphi(\alpha_{t^{-1}}(g))-\mu_\varphi(g)|\\&
\leq|\mu_\varphi(\alpha_{t^{-1}}(g))-(\varphi_*\zeta)(\alpha_{t^{-1}}(g))|
+|(\varphi_*\zeta)(g)-\mu_{\varphi}(g)|+\\&|(\varphi_*\zeta)(\alpha_{t^{-1}}(g))-((\varphi\circ
\sigma_t)_*\zeta)(g)|\\&\leq\frac{1}{n}.
\end{split}\end{equation*}

Now for any map $\sigma:G\to{\rm Sym}(d)$ with some $d\in
\mathbb{N},$ for $\xi\in D,$ define
\begin{equation*}
W_{\xi}(L_n,F_n,\delta_n,\sigma)=\left\{\varphi\in
Map(F_n,\delta_n,\sigma):\mu_{\varphi}=\xi\right\}.
\end{equation*}
Then

\noindent (i) $W_{\xi}(L_n,F_n,\delta_n,\sigma)\subset
Map_{\xi}(L_n,F_n,\frac{1}{3n},\sigma),$

\noindent (ii) $\bigcup\limits_{\xi\in
D}W_{\xi}(L_n,F_n,\delta_n,\sigma)=Map(F_n,\delta_{n},\sigma).$\\

Enumerate the measures in $D$ as follows:
\begin{equation*}
D=\{\xi_1,\xi_2,\cdots,\xi_{|D|}\}.
\end{equation*}
For any $\epsilon>0,$ there exists $\beta_i=\beta(\xi_i,\epsilon)\in
C^o_{W_{\xi_i}(L_n, F_n, \delta_n, \sigma)}$ and $
\beta_i\succcurlyeq\mathcal{U}^d$ such that

\begin{eqnarray*}\begin{split}
&\inf \{\sum\limits_{U\in\beta}\sup\limits_{\varphi\in U\cap
W_{\xi_i}(L_n,F_n,\delta_n,\sigma)}\exp(\sum\limits_{a=1}^{d}f(\varphi(a))):\beta\in
C^o_{W_{\xi_i}(L_n, F_n, \delta_n, \sigma)}{\rm
and~}\beta\succcurlyeq\mathcal{U}^d\}\\&\geq
\sum\limits_{U\in\beta_i}\sup\limits_{\varphi\in U\cap
W_{\xi_i}(L_n,F_n,\delta_n,\sigma)}\exp(\sum\limits_{a=1}^{d}f(\varphi(a)))e^{-\epsilon}.
\end{split}\end{eqnarray*}
Simplify $W_{\xi_i}(L_n,F_n,\delta_n,\sigma)$ as $W_i.$

Set
\begin{equation*}
\beta_i^*=\left\{V\cap
(X^d\backslash\bigcup\limits_{k:\overline{W}_k\cap\overline{W}_i=\emptyset}\overline{W}_k):V\in\beta_i\right\}
\end{equation*}

 Then $\beta_i^*\in C_{W _i}^o$ and
\begin{equation*}
\sum\limits_{U\in\beta_i}\sup\limits_{\varphi\in U\cap
W_{i}}\exp(\sum\limits_{a=1}^{d}f(\varphi(a)))\geq
\sum\limits_{U\in\beta^*_i}\sup\limits_{\varphi\in U\cap
W_{i}}\exp(\sum\limits_{a=1}^{d}f(\varphi(a))).
\end{equation*}

For any $V^{i*}\in \beta_i^*,$ if $\sup\limits_{\varphi\in
Map(F_n,\delta_n,\sigma)\cap
V^{i*}}\sum\limits_{a=1}^df(\varphi(a))>\sup\limits_{\varphi\in
W_k\cap V^{i*}}\sum\limits_{a=1}^df(\varphi(a)), $ then there exists
$k\in\{1,2,\cdots,|D|\}$ such that $\sup\limits_{\varphi\in
Map(F_n,\delta_n,\sigma)\cap
V^{i*}}\sum\limits_{a=1}^df(\varphi(a))=\sup\limits_{\varphi\in
W_k\cap V^{i*}}\sum\limits_{a=1}^df(\varphi(a)) $ with
$\overline{W_k}\cap\overline{W_i}\neq\emptyset.$ Choose
$\varphi^*\in\overline{W}_k\cap\overline{W}_i,$ there exist two
sequences $\{\varphi_j^k\}\subset W_k$ and  $\{\varphi_j^i\}\subset
W_i$ such that $\varphi_j^k\to \varphi^*, \varphi_j^i\to \varphi^*$
as $j\to\infty.$ Choose sufficiently large $j$ such that
$\rho_\infty(\varphi_j^k,\varphi^*)<\sqrt{\delta_n}$ and
$\rho_\infty(\varphi_j^i,\varphi^*)<\sqrt{\delta_n}.$ By the
definition of supremum, there exist $\varphi_1\in V^{i*}\cap W_k$
and $\varphi_2\in V^{i*}\cap W_i$ such that

\noindent (i) $|\sup\limits_{\varphi\in V^{i*}\cap
W_k}\sum\limits_{a=1}^df(\varphi(a))-\sum\limits_{a=1}^{d}f(\varphi_1(a))|<\frac{d}{9n},$
and

\noindent (ii) $|\sup\limits_{\varphi\in V^{i*}\cap
W_i}\sum\limits_{a=1}^df(\varphi(a))-\sum\limits_{a=1}^{d}f(\varphi_2(a))|<\frac{d}{9n}.$

Since $\varphi_1,\varphi_j^k\in W_k$ and $\varphi_2,\varphi_j^i\in
W_i,$ the following hold:

\noindent (iii)
$|\sum\limits_{a=1}^df(\varphi_1(a))-d\xi_k(f)|<d\delta_n,$

\noindent (iv)
$|\sum\limits_{a=1}^df(\varphi_j^k(a))-d\xi_k(f)|<d\delta_n,$

\noindent (v)
$|\sum\limits_{a=1}^df(\varphi_2(a))-d\xi_i(f)|<d\delta_n,$ and

\noindent (vi)
$|\sum\limits_{a=1}^df(\varphi_j^i(a))-d\xi_i(f)|<d\delta_n.$\\
Since $\rho_\infty(\varphi_j^k,\varphi^*)<\sqrt{\delta_n}$ and
$\rho_\infty(\varphi_j^i,\varphi^*)<\sqrt{\delta_n},$ we have

\noindent (vii)
$|\sum\limits_{a=1}^df(\varphi_j^i(a))-\sum\limits_{a=1}^df(\varphi^*(a))|<\frac{d}{6n},$
and

\noindent (viii)
$|\sum\limits_{a=1}^df(\varphi_j^k(a))-\sum\limits_{a=1}^df(\varphi^*(a))|<\frac{d}{6n}.$\\
Combining (i)-(viii) and trigonometrical inequalities, we have

\begin{equation*}
|\sup\limits_{\varphi\in Map(F_n,\delta_n,\sigma)\cap
V^{i*}}\sum\limits_{a=1}^df(\varphi(a))-\sup\limits_{\varphi\in
W_i\cap V^{i*}}\sum\limits_{a=1}^df(\varphi(a))|<\frac{d}{n}.
\end{equation*}

Hence,

\begin{eqnarray*}\begin{split}
&N^*(\sigma,F_n,\delta_n,\mathcal{U},f)\\&
\leq\sum\limits_{i=1}^{|D|}\sum\limits_{U\in
\beta_i^*}\sup\limits_{\varphi\in U\cap
W_i}\exp(\sum\limits_{a=1}^{d}f(\varphi(a))+\frac{d}{n}).
\end{split}\end{eqnarray*}
There exists $\nu=\nu_\sigma\in D$ such that

\begin{eqnarray*}\begin{split}
&N^*(\sigma,F_n,\delta_n,\mathcal{U},f)\\&\leq |D|\sum\limits_{U\in
\beta(\nu,\epsilon)}\sup\limits_{\varphi\in U\cap
W_\nu(L_n,F_n,\delta_n,\sigma)}\exp(\sum\limits_{a=1}^{d}f(\varphi(a))+\frac{d}{n})\\&\leq|D|\inf
\{\sum\limits_{U\in\beta}\sup\limits_{\varphi\in U\cap
W_\nu(L_n,F_n,\delta_n,\sigma)}\exp(\sum\limits_{a=1}^{d}f(\varphi(a))+\frac{d}{n}):\beta\in
C^o_{W_\nu(L_n,F_n,\delta_n,\sigma)}{\rm~
and~}\beta\succcurlyeq\mathcal{U}^d\}e^\epsilon\\& \leq
|D|\exp(\nu(f)d+\delta_n d+\frac{d}{n})\inf
\{\sum\limits_{U\in\beta}1: \beta\in
C^o_{W_\nu(L_n,F_n,\delta_n,\sigma)}{\rm~
and~}\beta\succcurlyeq\mathcal{U}^d\}e^\epsilon
\\&\leq
|D|\exp(\nu(f)d+\delta_n d+\frac{d}{n})N(\mathcal{U}^d,
Map_\nu(L_n,F_n,\frac{1}{3n},\sigma))e^\epsilon.
\end{split}\end{eqnarray*}
Letting $\sigma$ run through the terms of the sofic approximation
sequence $\Sigma,$ by the pigeonhole principle there exists
$\mu_n\in D$ and a sequence $i_1<i_2<\cdots$ such that the following
hold:

\noindent (i) $
P^*(F_n,\delta_n,\mathcal{U},f)=\lim\limits_{k\to\infty}\frac{1}{d_{i_k}}\log
N^*(\sigma_{i_k},F_n,\delta_n,\mathcal{U},f); $

\noindent (ii)
\begin{equation*}\begin{split}
&\frac{1}{d_{i_k}}\log
N^*(\sigma_{i_k},F_n,\delta_n,\mathcal{U},f)\\&\leq\frac{1}{d_{i_k}}\log(|D|\exp(\mu_n(f)d_{i_k}+\delta_n
d_{i_k}+\frac{d_{i_k}}{n})N(\mathcal{U}^{d_{i_k}},
Map_{\mu_n}(L_n,F_n,\frac{1}{3n},\sigma_{i_k}))e^\epsilon),
\end{split}\end{equation*}
for all $k\in \mathbb{N};$

\noindent (iii) $\max\limits_{t\in F_n, g\in
L_n}|\mu_n(\alpha_{t^{-1}}(g))-\mu_n(g)|<\frac{1}{n}.$

\noindent Then $h_{\mu_n}(L_n, F_n,
\frac{1}{3n},\mathcal{U})+\mu_n(f)+\frac{10}{9n}\geq
P^*(F_n,\delta_n,\mathcal{U},f).$

Let $\mu$ be a weak* limit point of the sequence
$\{\mu_n\}_{n\in\mathbb{N}}.$ If $t\in G$ and $g\in \{g_m\}_{m\in
\mathbb{N}}, $ then
\begin{equation*}\begin{split}
&|\mu(\alpha_{t^{-1}}(g))-\mu(g)|\\&\leq
|\mu(\alpha_{t^{-1}}(g)))-\mu_n(\alpha_{t^{-1}}(g))|+|\mu_n(\alpha_{t^{-1}}(g))-\mu_n(g)|+|\mu_n(g)-\mu(g)|.
\end{split}\end{equation*}
Since the infimum of the right hand side over all $n\in\mathbb{N}$
is 0 and $\{g_m\}_{m\in\mathbb{N}}$ is dense in $C(X),$ we deduce
that $\mu$ is $G$-invariant.

Let $F\in \mathcal{F}(G), L\in \mathcal{F}(C(X))$ and $\delta>0.$
Take an integer $n$ such that $F\subset F_n, \frac{1}{3n}\leq
\frac{\delta}{4}, \max\limits_{g\in
L\cup\{f\}}|\mu_n(g)-\mu(g)|<\frac{\delta}{6}$ and for any $g\in L,$
there exists $g'\in L_n$ such that $\|g-g'\|_\infty\leq
\frac{\delta}{4}.$ Then for any map $\sigma$ from $G$ to Sym$(d)$
for some $d\in\mathbb{N}, \varphi\in Map_{\mu_n}(L_n, F_n,
\frac{1}{3n},\sigma)$ and $g\in L,$ we have

\begin{equation*}\begin{split}
&|(\varphi_*\zeta)(g)-\mu(g)|\\&\leq|(\varphi_*\zeta)(g)-(\varphi_*\zeta)(g')|+|(\varphi_*\zeta)(g')-\mu_n(g')|+|\mu_n(g')-\mu_n(g)|+|\mu_n(g)-\mu(g)|\\&
<\frac{3\delta}{4}+\frac{1}{3n}\leq \delta,
\end{split}\end{equation*}
and hence $\varphi\in Map_\mu(L, F, \delta, \sigma).$ Thus
$Map_{\mu_n}(L_n, F_n, \frac{1}{3n},\sigma)\subset Map_{\mu}(L, F,
\delta, \sigma)$ and then

\begin{equation*}\begin{split}
&h_{\mu}(L, F, \delta, \mathcal{U})+\mu(f)\\&\geq h_{\mu_n}(L_n,
F_n, \frac{1}{3n},\mathcal{U})+\mu_n(f)-\frac{\delta}{6}\\&\geq
P^*(F_n,\delta_n,\mathcal{U},f)-\frac{10}{9n}-\frac{\delta}{6}\\&\geq
P^*(\mathcal{U},f)-\delta.
\end{split}\end{equation*}
Consequently, $P^*(\mathcal{U},f)\leq h_{\mu}(\mathcal{U})+\mu(f).$

\begin{lem}
If $\mathcal{U}\in C_X^o, f\in C(X),$ then
$P(\mathcal{U},f)\geq\max\limits_{\mu\in
M(X,G)}\{h_\mu(\mathcal{U})+\mu(f)\}.$
\end{lem}
\noindent{\it Proof.} Let $F\in \mathcal{F}(G),$ and $\delta>0.$ Set
$L_1=\{f\}.$  There exists $\beta\in
C_{Map(F,\frac{\delta}{2},\sigma_i)}$and
$\beta\succcurlyeq\mathcal{U}^{d_i}$ such that

\begin{eqnarray*}\begin{split}
&N(\sigma_i,F,\frac{\delta}{2},\mathcal{U},f)\exp(\frac{\delta}{2})\\&\geq
\sum\limits_{U\in\beta}\sup\limits_{\varphi\in U\cap
Map(F,\frac{\delta}{2},\sigma_i)}\exp(\sum\limits_{a=1}^{d_i}f(\varphi(a)))\\&\geq
\sum\limits_{U\in\beta}\sup\limits_{\varphi\in U\cap
Map_\mu(L_1,F,\frac{\delta}{2},\sigma_i)}\exp(\sum\limits_{a=1}^{d_i}f(\varphi(a)))\\&
\geq
N(\mathcal{U}^{d_i},Map_\mu(L_1,F,\frac{\delta}{2},\sigma_i))\exp(d_i
\mu(f)-d_i\frac{\delta}{2}).
\end{split}\end{eqnarray*}
Thus,
\begin{equation*}
P(F, \frac{\delta}{2},\mathcal{U},f) \geq
h_{\mu}(L_1,F,\frac{\delta}{2},\mathcal{U})+\mu(
f)-\frac{\delta}{2}.
\end{equation*}

\begin{equation*}
P(F, \frac{\delta}{2},\mathcal{U},f) +\frac{\delta}{2}\geq
h_{\mu}(\mathcal{U})+\mu( f).
\end{equation*}

Letting $\delta\to0,$  we have
\begin{equation*}
P(\mathcal{U},f)\geq h_{\mu}(\mathcal{U})+\mu(f).
\end{equation*}

\begin{thm}\label{thm3.1}{\rm (variational priciple)}
If $\mathcal{U}\in C_X^o, f\in C(X),$ then
$P(\mathcal{U},f)=P^*(\mathcal{U},f)=\max\limits_{\mu\in
M(X,G)}\{h_\mu(\mathcal{U})+\mu(f)\}.$
\end{thm}

The following corollary can be obtained directly from the
variational principle.

\begin{cro}\label{cor3.1}
If $\mathcal{U}\in C_X^o, \mu\in M(X,G),$ then
$h_\mu(\mathcal{U})\leq\inf\limits_{f\in
C(X)}\left\{P(\mathcal{U},f)-\mu(f)\right\}.$
\end{cro}
\section{Properties of local sofic pressure}

In this section, we are to compare local sofic pressure with
classical counterparts in the setting of the group being amenable
and to study properties of local sofic pressure. This section makes
full use of variational principles.

A group $G$ is said to be amenable if for all $\epsilon>0$ and all
$K\in\mathcal{F}(G),$ there exists $F\in \mathcal{F}(G) $ such that
\begin{equation*}
\frac{|F\Delta KF|}{|F|}<\epsilon.
\end{equation*}
A countable group is amenable if and only if there is a sequence
$\{F_n\}_{n\in \mathbb{N}}\subset\mathcal{F}(G)$ such that
\begin{equation*}
\lim\limits_{n\to\infty}\frac{|F_n\Delta KF_n|}{|F_n|}=0,
\end{equation*}
for all $K\in\mathcal{F}(G).$ Such a sequence is called a
F$\phi$lner sequence of $G.$

Now, we give the definition of local amenable pressure which is a
special case in \cite{[LiaYan]}.

Let $(X,G)$ be a system, where $G$ is a countable infinite amenable
group. For $\mathcal{U}, \mathcal{V}\in C_X,$ set
$\mathcal{U}\bigvee\mathcal{V}=\{U\cap V:U\in
\mathcal{U},V\in\mathcal{V}\}.$ Given $F\in\mathcal{F}(G)$ and
$\mathcal{U}\in C_X,$ set $\mathcal{U}_F=\bigvee_{g\in
F}g^{-1}\mathcal{U}.$

Let $\nu\in M(X)$ and $ \alpha\in \mathcal{P}_X.$ Set
$H_\nu(\alpha)=-\sum\limits_{A\in\alpha}\nu(A)\log\nu(A),$ and
$H_\nu(\mathcal{V})=\inf\limits_{\alpha\in
\mathcal{P}_X,\alpha\succcurlyeq\mathcal{V}}H_\nu(\alpha).$

For $E\in \mathcal{F}(G),\mathcal{U}\in C_X^o$ and $f\in C(X),$
define
\begin{equation*}
P^a_E(\mathcal{U},f)=\inf\left\{\sum\limits_{V\in\mathcal{V}}\sup\limits_{x\in
V}e^{\sum\limits_{g\in E}f\circ g(x)}:\mathcal{V}\in
C_X{\rm~and~}\mathcal{V}\succcurlyeq\mathcal{U}_E\right\}.
\end{equation*}
The existence of the following limits are due to Ornstein and Weiss
\cite{[OrnWei]}.
\begin{equation*}
P^a(\mathcal{U},f)=\lim\limits_{n\to\infty}\frac{1}{|F_n|}\log
P^a_{F_n}(\mathcal{U},f),
\end{equation*}
where $\{F_n\}_{n\in\mathbb{N}}$ is a F$\phi$lner sequence.

\noindent Let $\mathcal{U}\in C_X^o$ and $\mu\in M(X,G).$ Set
\begin{equation*}
h_\mu^a(\mathcal{U})=\lim\limits_{n\to\infty}\frac{1}{|F_n|}H_\mu(\mathcal{U}_{F_n})
\end{equation*}
where $\{F_n\}_{n\in\mathbb{N}}$ is a F$\phi$lner sequence.

\begin{thm}
If $G$ is a countable amenable group. Let $\mathcal{U}\in C_X^{o},
f\in C(X),$ then $P(\mathcal{U},f)=P^{a}(\mathcal{U},f).$
\end{thm}
\noindent{\it Proof.}
\begin{equation*}\begin{split}
P(\mathcal{U},f)&=\sup\limits_{\mu\in
M(X,G)}\left\{h_{\mu}(\mathcal{U})+\mu(f)\right\}\\&
=\sup\limits_{\mu\in
M(X,G)}\left\{h^a_{\mu}(\mathcal{U})+\mu(f)\right\}\\&
=P^a(\mathcal{U},f).
\end{split}\end{equation*}
The second equality depends on the fact that
$h^a_{\mu}(\mathcal{U})=h_{\mu}(\mathcal{U})$ (see \cite{[Zha]}).
The third equality relies on the variational principle in
\cite{[LiaYan]}.

\begin{prop}
Let $\mathcal{U}\in C_X^o.$ If $f, g\in C(X), s\in G$ and $c\in
\mathbb{R},$ then the following hold:

{\rm (i)} $P(\mathcal{U}, f+c)=P(\mathcal{U}, f)+c,$

{\rm (ii)} $P(\mathcal{U}, f+g)\leq P(\mathcal{U}, f)+P(\mathcal{U},
g),$

{\rm (iii)} If $f\leq g,$ then $P(\mathcal{U},f )\leq
P(\mathcal{U},g ),$ in particular, $h(G,\mathcal{U})+\min f\leq
P(\mathcal{U},f)\leq h(G,\mathcal{U})+\max f,$ and $ P(\mathcal{U},
f)\leq P(\mathcal{U}, |f|),$

{\rm (iv)} $P(\mathcal{U},\cdot)$ is either finite valued or
constantly $\pm\infty,$

{\rm (v)} If $P(\mathcal{U},\cdot)\neq\pm\infty,$ then
$\|P(\mathcal{U}, f)-P(\mathcal{U}, g)\|\leq\|f-g\|_\infty,$

{\rm (vi)} $P(\mathcal{U}, f+g\circ\alpha_s-g)=P(\mathcal{U},f),$

{\rm (vii)} Assume that $P(\mathcal{U},\cdot)\neq-\infty,
P(\mathcal{U}, cf)\leq c\cdot P(\mathcal{U}, f)$ if $c\geq1$ and
$P(\mathcal{U}, cf)\geq c\cdot P(\mathcal{U}, f)$ if $c\leq1.$

{\rm (viii)} $P(\mathcal{U},\cdot)$ is convex.
\end{prop}

\noindent{\it Proof.} (i), (ii), (iii), (vi), (vii) and (viii) can
be obtained from variational principle (Theorem \ref{thm3.1})
easily. (iv) are from (iii), and (v) are from (i) and (iii).

\noindent Recall that a finite signed measure on $X$ is a map
$\mu:\mathcal{B}(X)\to\mathbb{R}$ which is countably additive, where
$\mathcal{B}(X)$ is the $\sigma$-algebra of Borel subsets of $X.$
(See \cite{Wal2}P$_{221}$)
\begin{thm}\label{thm4.2}
Assume that $h(G,\mathcal{U})\neq\pm\infty.$ Let
$\mu:\mathcal{B}(X)\to\mathbb{R}$ be a finite signed measure. If
$\mu(f)\leq P(\mathcal{U}, f)$ for all $f\in C(X),$ then $\mu\in
M(X,G).$
\end{thm}
\noindent{\it Proof.} Firstly, show $\mu$ takes only non-negative
values. Suppose $f\geq0.$ If $\kappa>0$ and $n>0$ we have

\begin{equation*}\begin{split}
&\mu(n(f+\kappa))\\& =-\mu(-n(f+\kappa))\geq
-P(\mathcal{U},-n(f+\kappa))\\&\geq-(h(G,\mathcal{U})+\max(-n(f+\kappa)))\\&=-h(G,\mathcal{U})
+n\min(f+\kappa)\\&>0{\rm~ for~ large~}n.
\end{split}\end{equation*}
Therefore $\mu(f+\kappa)>0$ and hence $\mu(f)\geq0$ as desired.

Secondly, show $\mu(X)=1.$ If $n\in\mathbb{Z}$ then $\mu(n)\leq
P(\mathcal{U},n)= h(G,\mathcal{U})+n,$ so that $\mu(X)\leq1+
\frac{h(G,\mathcal{U})}{n}$ if $n>0$ and hence $\mu(X)\leq1,$ and
$\mu(X)\geq1+ \frac{h(G,\mathcal{U})}{n}$ if $n<0$ and hence
$\mu(X)\geq1.$ Therefore $\mu(X)=1.$

Thirdly,  we show $\mu\in M(X,G).$ Let $s\in G, n\in\mathbb{Z}$ and
$f\in C(X). n\mu(f\circ\alpha_s-f)\leq
P(\mathcal{U},n(f\circ\alpha_s-f))=h(G,\mathcal{U}).$ If $n>0$ then
dividing both sides by $n$ and letting $n$ go to $\infty$ yields
$\mu(f\circ\alpha_s-f)\leq0,$ and if $n<0$ then dividing both sides
by $n$ and letting $n$ go to $-\infty$ yields
$\mu(f\circ\alpha_s-f)\geq0.$ Therefore
$\mu(f\circ\alpha_s)=\mu(f),$ for any $f\in C(X),s\in G.$ Thus
$\mu\in M(X,G).$

Now, we end this paper with the following questions:

(i) Let $\mathcal{U}\in C_X^o,$ and $\mu\in M(X,G).$ Let
$\mu=\int_{M^e(X,G)}m d\tau(m)$ be the ergodic decomposition of
$\mu.$ Do we have $h_\mu(\mathcal{U})=\int_{M^e(X,G)}
h_{m}(\mathcal{U})d\tau(m)?$

(ii) Does  Theorem \ref{thm4.2} still hold without the hypothesis
$h(G,\mathcal{U})\neq\pm\infty?$

(iii) Does the opposite inequality of Corollary \ref{cor3.1} hold?

In fact, If (i) is true, then it is easy to prove that (iii) holds.


\noindent {\bf Acknowledgements.}  We would like to thank Prof.
Hanfeng Li for some comments. The
work was supported by the National Natural Science Foundation of
China (10971100) and National Basic Research Program of China (973
Program) (2007CB814800).


\begin{thebibliography}{50}
\bibitem{[Bow2]} L. Bowen, Measure conjugacy invariants for actions of countable sofic groups,
J. Amer. Math. Soc., {\bf 23} (2010), no. 1, 217-245.

\bibitem{[Bow1]} L. Bowen, Sofic entropy and amenable groups, Ergod.
Th. Dynam. Sys., to appear.

\bibitem{[Bow0]} R. Bowen, Topological entropy for non-compact set,
Trans. Amer. Math. Soc., {\bf 184} (1973), 125-136.


\bibitem{[Chu]} N. Chung, The variational principle of
topological pressure for actions of sofic groups. arXiv:
1110.0699v1.

\bibitem{[GlaWei]} E. Glasner \& B. Weiss, On the interplay between
measurable and topological dynamics, in Handbook of dynamical
systems, Vol. 1B, Elsevier B. V., Amsterdam, (2006) 597-648.

\bibitem{[Goo1]} T. Goodman, Topological entropy bounds
measure-theoretic entropy, Proc. Amer. Math. Soc., {\bf 23} (1969),
679-688.

\bibitem{[Goo2]} L. Goodwyn, Relating topological entropy and
measure entropy, Bull. London Math. Soc., {\bf 3} (1971), 176-180.

\bibitem{[HuaYe]} W. Huang \& X. Ye, A local variational relation
and applications, Israel J. Math., {\bf 151} (2006), 237-279.

\bibitem{[HuaYeZha]} W. Huang, X. Ye \& G. Zhang, Local entropy
theory for a countable discrete amenable group action, J. Funct.
Anal., {\bf 261}(2011), 1028-1082.

\bibitem{[HuaYi]} W. Huang \& Y. Yi, A Local variational principle
of pressure and its applications to equilibrium states, Israel J.
Math., {\bf 161} (2007), 29-74.

\bibitem{[KerLi1]} D. Kerr \& H. Li, Entropy and the variational
principle for actions of sofic groups, Invent. Math., {\bf 186}
(2011),  501--558.

\bibitem{[KerLi2]} D. Kerr \& H. Li, Soficity, amenability, and
dynamical entropy, Amer. J. Math., to appear.

\bibitem{[Li]} H. Li, Sofic mean dimension, arXiv: 1105.0140v1.


\bibitem{[LiaYan]} B. Liang \& K. Yan, Topological pressure
for sub-additive potentials of amenable group actions, J. Func.
Anal., {\bf 262} (2012), 584-601.


\bibitem{[Mis]} M. Misiurewicz, A short proof of the variational
principle for a $\mathbb{Z}^n_+$ action on a compact space, Bull.
Acad. Polon. Sci. Sr. Sci. Math. Astronom. Phys.,  {\bf 24} (1976),
1069-1075.

\bibitem{[OrnWei]} D. Ornstein \& B. Weiss, Entropy and isomorphism
theorems for actions of amenable groups, J. Anal. Math., {\bf 48}
(1987), 1-141.


\bibitem{[PesPit]}Y. Pesin \& B. Pitskel,  Topological pressure and
the variational principle for noncompact sets,  Func. Anal. Appl.,
 {\bf18} (1984), 307-318.


\bibitem{[Rom]} P. Romagnoli,  A local variational principle for the topological entropy,
Ergod. Th. Dynam. Sys., {\bf 23} (2003), 1601-1610.

\bibitem{[Wal1]}P. Walter, A variational principle for the pressure
of continuous transformations, Amer. J. Math., {\bf 97} (1975),
937-971.


\bibitem{Wal2} P. Walters,  An Introduction to Ergodic Theory, Springer, New York, 1982.



\bibitem{[Zha]} G. Zhang, Local variational principle concering entropy of a
sofic group action, J. Func. Anal., to appear.















\end{thebibliography}
\end{document}